\documentclass[11pt]{article}
\usepackage{authblk}
\usepackage[numbers]{natbib}
\usepackage{preamble}
\usepackage[colorlinks=true, allcolors=blue]{hyperref}
\usepackage[
    left=3cm,   
    right=3cm,  
    top=2.25cm,    
    bottom=2.25cm, 
]{geometry}

\author[1]{Charly Robinson La Rocca\thanks{Email: charly.robinson.la.rocca@umontreal.ca; Corresponding author.}}
\author[2]{Jean-François Cordeau}
\author[1]{Emma Frejinger}

\affil[1]{Department of Computer Science and Operations Research, Université de Montréal, Canada}

\affil[2]{Department of Logistics and Operations Management, HEC Montréal, Canada}

\begin{document}

\date{\today}
\title{One-shot Learning for MIPs with SOS1 Constraints}
\maketitle

\begin{abstract}
    Efficient algorithms and solvers are required to provide optimal or near-optimal solutions quickly and enable organizations to react promptly to dynamic situations such as supply chain disruptions or changing customer demands. State-of-the-art mixed-integer programming (MIP) solvers are crafted to tackle a wide variety of problems, yet many real-world situations are characterized by problem instances that originate from a narrow distribution. This has inspired the creation of tailored approaches that exploit historical data to inform heuristic design. Deep learning (DL) methods are typically used in this context to extract patterns from data, but they require large datasets and comprehensive hyperparameter tuning for strong performance. This article describes a one-shot learning heuristic that leverages solutions discovered within the branch-and-bound tree to construct a model with minimal overhead. We evaluate our method on the locomotive assignment problem (LAP) and sets of MIPLIB instances that contain constraints based on special ordered sets of type 1. Experimental results include a comparison with multiple primal heuristics and state-of-the-art MIP solvers. We show that the method is most effective with CPLEX in terms of the average primal gap.
\end{abstract}


\renewcommand{\V}{\mathcal{V}}
\newcommand{\Vr}{\mathcal{V}'}
\renewcommand{\P}{\mathcal{P}}
\renewcommand{\Pr}{\mathcal{P}'}
\renewcommand{\K}{\mathcal{K}}
\newcommand{\kt}{k'}
\newcommand{\ko}{k^*}
\newcommand{\x}{\mathbf{x}}
\newcommand{\ct}{c^T}
\newcommand{\xv}{x_v}
\newcommand{\xvk}{x_v^k}
\newcommand{\xlp}{\tilde{\mathbf{x}}}
\newcommand{\xvklp}{\tilde{x}_v^k}

\newcommand{\xvkt}{x_{vt}^k}
\newcommand{\kvt}{k_{vt}}
\newcommand{\kp}{\mathbf{k}_{v}}
\newcommand{\X}{\mathcal{X}_p}
\newcommand{\tp}{T_p}
\renewcommand{\r}{r}
\newcommand{\bnb}{B\&B}
\newcommand{\xc}{\mathbf{x}^*}
\newcommand{\xinc}{\mathbf{x}_{\text{inc}}}
\newcommand{\dmax}{d_{\text{max}}}

\newcommand{\cplex}{\text{CPLEX}}
\newcommand{\pnf}{\text{PNF} }
\newcommand{\runtime}{\text{runtime}}
\newcommand{\solve}{\text{solve}}
\newcommand{\probe}{\text{probe}}
\newcommand{\select}{\text{select}}
\newcommand{\score}{\text{score}}
\newcommand{\predict}{\text{predict}}

\section{Introduction}
The intersection of mixed-integer programming (MIP) and machine learning (ML) presents a unique opportunity to advance the frontiers of optimization and tackle complex problems in areas such as transportation~\cite{sun_generalization_2021}, energy~\cite{xavier_learning_2021}, and finance~\cite{fischer_deep_2018}. This trend has emerged in recent years with the maturation of specialized hardware \cite{wang2019benchmarking} and software \cite{paszke2019pytorch} stacks. Combined with the availability of large datasets,   ML can exploit patterns from data to design tailored algorithms. Today, there are many strategies that take advantage of learning methods to improve upon standard solvers.  Bengio et al. \cite{Bengio_2020} provide a tour d'horizon on ML for combinatorial optimization (CO) and describe  the different ways they can be integrated.  The first idea uses ML to directly predict solutions of the problem instances from the input features. One example of this is the application of the pointer network \cite{vinyals2015pointer} to learn to predict an optimal permutation for the travelling salesman problem (TSP). The second approach is to learn to optimize the hyperparameters of a CO algorithm \cite{xu2011hydra}. Third, ML can also be used alongside an optimization algorithm. In practice, this is often done via the integration of learning in the branch-and-bound (B\&B) framework \cite{lodi2017learning}. Our approach belongs to the first class of hybrid methods.

This paper proposes a primal heuristic for MIPs with special ordered sets of type 1 (SOS1 or S1) \cite{beale1970special} that we refer to as Probe and Freeze\footnote{Source code available at https://github.com/laroccacharly/MLSOS} (PNF). \startblue Our objective is to find high-quality solutions for instances that are both recurring and challenging to solve for general-purpose solvers. We aim to craft a heuristic that is easy to reproduce so it can be used as a baseline for future research. Our experimental setup specifically considers the effect of the problem reduction as a preprocessing step instead of integrating it into the branch-and-bound. This decision was made to have a solver-agnostic approach and simplify the algorithm design. 
\stopblue 

At a high level, S1 constraints are modelled using a sum over a set of binary variables constrained to be 1. We use the shorter notation S1 to refer to special ordered sets of type 1 in the remainder of this paper.  The idea of PNF is to use the probing data to train a classifier to predict which variable should take the value 1 in each S1 constraint of the MIP.  Training is done in a one-shot fashion, analogous to the one-shot learning (OSL) paradigm, which is the main difference compared to other works in the literature. The method uses entropy, a key concept in information theory, to quantify uncertainty. The predictions are used to freeze a subset of the variables in the S1 constraints to produce a reduced problem that can be significantly easier to solve.  Several MIP heuristics adopt the strategy of fixing variables and solving the reduced problem \cite{danna2005exploring, berthold_rens_2014, wolsey2020integer}. Our heuristic can be seen as a form of Large Neighborhood Search (LNS) \cite{shaw_using_1998} where the neighborhood is defined by the variables that are not fixed.  

\textbf{Motivation for S1.} The motivation for this approach comes from the hypothesis that S1 constraints have a structuring effect on the final solution.  These constraints are useful when a choice involves multiple options or resources, and only one can be selected. For example, in the facility location problem (FLP), the objective is to determine the optimal locations for a set of facilities to serve a given set of customers. The S1 constraints can either be used to model which type of facility to build or how to assign customers to facilities. In the locomotive assignment problem (LAP), which is a problem studied in this paper, S1 constraints model the assignment of consists (i.e., groups of locomotives) to trains. The flow variables in the LAP are often trivial to optimize for a given assignment of consists to trains. In other words, the S1 constraints are the main source of complexity in the problem. Furthermore, the S1 constraint can be modified to model more complex types of assignments. \startblue  For example, the S2 constraint is a possible extension where two adjacent variables can take non-zero values. \stopblue 

\textbf{Motivation for OSL.} The main motivation behind OSL is the cost associated with the data generation process. In the context of MIPs, the training dataset typically comes from solving a large number of instances. This is a time-consuming process that requires a significant amount of computational resources.  Another reason is the fact that most of the rich data structure generated by B\&B is usually overlooked during the search. This work explores how to make use of this readily available data using statistical tools to better inform a heuristic. As a side effect of using OSL, we can produce a method that is easy to understand and reproduce. It is also worth noting that the OSL paradigm is not a hard requirement for our approach. In fact, we can use any ML classifier to determine which variable the value 1 should be assigned to. The OSL model acts as a strong baseline which can be combined with other modern ML models to improve the accuracy of predictions. 

\textbf{Motivation for our test sets.} The LAP is a problem that is not only challenging to solve but also recurring, i.e., similar problem instances are solved repeatedly over time. Given the large scale of the problem, there is a substantial economic incentive to find a heuristic that can solve the LAP quickly.  We use the Canadian National Railway Company (CN) as the subject of our case study for the LAP. We test our method on two sets of instances with different difficulty levels to assess the capabilities of our approach. Furthermore, we include a comparison with a set of MIPLIB instances with S1 constraints to evaluate the performance of our method on a more diverse set of problems. 

This document first presents some related work on the integration of ML in MIP solvers. This is followed by an overview of our methodology which is composed of three routines: probe, select and freeze. Finally, we include results for the primal gap and runtime on sets of LAP and MIPLIB instances. For each instance, we compare the performance of PNF with different primal heuristics and state-of-the-art MIP solvers: CPLEX \cite{cplex}, Gurobi \cite{gurobi} and SCIP \cite{scip}. Both CPLEX and Gurobi are commercial solvers, while SCIP is an open-source solver. Commercial solvers are often considered the gold standard for MIPs, but they are not always available to researchers and, when they are, it is not always easy to know exactly what procedures are applied by the solver during the solution process. Hence, we include SCIP in our comparison because it is commonly used in related works and it is more transparent than commercial solvers.

\section{Related Works}
 The underlying algorithm behind state-of-the-art (SOA) solvers is the ``divide and conquer’’  framework called branch-and-bound \cite{LandDoig:BranchAndBound}. Alone, the classical B\&B implementation often fails to close the optimality gap in an acceptable amount of time. Extensive research has been done to improve the performance of B\&B. In the following, we first present some of the most relevant works in the operational research (OR) literature. We then discuss recent advances in ML for MIPs. 

\textbf{OR approaches.} Heuristics for MIPs can be classified in two general classes: constructive and improvement heuristics \cite{hanafi2017mathematical}. An example of a constructive heuristic  is the Feasiblity Pump (FP) \cite{fischetti2005feasibility, bertacco2007feasibility, achterberg2007improving} where the main goal is to quickly produce feasible solutions. Typically, solutions generated by this type of heuristic do not come with guarantees. Improvement heuristics iteratively apply updating steps to an initial solution (feasible or not) to refine the quality of solutions.   In that category, there are LNS heuristics \cite{shaw_using_1998, ropke_adaptive_2006, pisinger_large_2019} that define a neighborhood of interest and explore that space to improve the current incumbent solution. The LNS framework has been studied extensively and there are many variants such as Relaxation Induced Neighborhood Search (RINS) \cite{danna2005exploring}, Relaxation Enforced Neighborhood Search (RENS) \cite{berthold_rens_2014}. The former defines a neighborhood related to the matching set between the LP relaxation and the incumbent solution. In the latter heuristic, the neighborhood focuses on the feasible roundings of the LP relaxation. The Local Branching (LB) heuristic \cite{fischetti_local_2003} introduced the so-called \emph{local branching} constraints to maintain proximity to the incumbent solution. Its distinctive feature is the utilization of a pseudo-cut for soft variable fixing.  Instead of adding a constraint, Proximity Search \cite{fischetti_proximity_2014} uses a penalty term in the objective function to minimize the distance to the incumbent solution. More recently, heuristics that employ ML have been developed to determine dynamically which LNS heuristic to run \cite{hendel2022adaptive}. 

Constructive heuristics, also called Start Heuristics (SH) \cite{berthold2006primal}, are most similar to the methodology proposed in this paper. Diving \cite{lazic2016variable} and Rounding \cite{achterberg2012rounding, neumann2019feasible} heuristics are two other examples that fit in that category. Diving heuristics employ a depth-first-search (DFS) strategy which focuses on decreasing the number of fractional variables in the LP. At every iteration, the algorithm selects a fractional variable to bound. Fractional, coefficient and pseudocost diving all share this idea but they use a different strategy to select which variable to bound \cite{berthold2006primal}. Rounding heuristics use the notion of up and down locks to perform a rounding that guarantees that all linear constraints are satisfied \cite{achterberg2012rounding}. The Zero Integer (ZI) round \cite{wallace2010zi} heuristic calculates the bounds for fractional variables such that the LP stays feasible. Variables are shifted to the corresponding bound that minimizes fractionality. 

Two relevant references  \cite{gilpin_information-theoretic_2011, karzan2009information} in the OR literature exploit information theory principles to aid decision-making during branching.   In this context, the LP value of a variable is interpreted as a probability.  Similar to our approach, the first paper \cite{gilpin_information-theoretic_2011} quantifies the amount of uncertainty of a binary variable using the notion of entropy \cite{shannon1948mathematical}. They describe Entropic Branching (EB), a look-ahead strategy similar to strong branching, where variable selection is done in a way that minimizes uncertainty in child nodes.  The second paper \cite{karzan2009information} uses a restart strategy to collect partial assignments from fathomed nodes called clauses. They experiment with both branching alternatives and clause inequalities that use this information to guide the search. This approach is quite similar to the Rapid Learning heuristic \cite{berthold2010rapid, berthold2019local} (related to no-good or conflict learning) which applies a depth-first search to collect valid conflict constraints. From this probing step, it identifies explanations for the infeasibility and exploits that information in the remainder of the search to prune parts of the tree.

\textbf{ML for MIPs.} In the previous paragraphs, we introduced heuristics that are readily available within SOA solvers. The solvers typically tune the schedules of these heuristics using the empirical performance on a broad set of test instances. The works of \cite{schedulingheuristics} use a data-driven approach to obtain a problem-specific schedule for different heuristics. By accomplishing this, they are able to improve the primal integral by 49\% on two classes of instances.  The scheduling of heuristics is just one of many possible ways to improve upon the default behaviour of SOA solvers. Other papers use ML to improve the node and variable selection strategies \cite{gasse2019exact, khalil}.  Learning to cut is another promising approach for strengthening the linear relaxation \cite{tang2020reinforcement}. 

Graph neural networks (GNN) \cite{hamilton2020graph} are commonly used to model the graph representation of a MIP. This architecture is popular because a MIP can be represented as a bipartite graph on variables
and constraints. They can be trained via supervised learning to predict the stability \cite{ding2020accelerating}, the bias \cite{khalil2022mip} or the value \cite{huang2022improving} of binary variables. The predictions are used to guide the search in the B\&B tree.  These approaches keep the exactness of the algorithm because all possible decisions made by the ML model are valid by design. The optimality certificate is often obtained by iteratively improving both the upper and lower bounds until they converge to the same solution.  However, for many practical applications, it is preferred to have a fast heuristic that can quickly produce a high-quality solution without the need to prove optimality. For example, the Deep Learning Tree Search (DLTS) \cite{hottung2020deep} algorithm uses neural network predictions for  both branching and bounding. Given that it relies on an estimated bound instead of a proven one, this approach is not guaranteed to converge to the optimal solution. Stochastic policies trained to generate solutions to CO problems have also been proposed \cite{toenshoff2021graph}. This type of approach iteratively generates solutions by sampling from a learned distribution over the variables, hence it does not produce any optimality certificate. Nonetheless, the authors show that the solutions generated by their approach are of high quality.

Our work differs from previous hybrid approaches in two key aspects. First, we fix a subset of the variables in the MIP formulation (before calling the solver) instead of directly guiding the search. In that sense, our approach is most similar to the Relax and Fix \cite{wolsey2020integer} primal heuristic where variables are fixed based on the fractional solution of the LP.  Second, we do not rely on a large dataset to train the model. Instead, analogous to one-shot learning, training is done on the probing data which contains few samples. 

One-shot learning has emerged as a significant direction in the ML literature, aiming to mimic human-like learning from few examples.  There is a vast literature on the subject and we refer to \cite{wang_generalizing_2021} for a related survey. We cover here two notable approaches that are well known because they are among the first to demonstrate the potential of OSL for classification problems.  Matching Networks \cite{vinyals2016matching} employ an attention mechanism over a learned embedding of the training set of examples to predict classes for unlabeled points, functioning as a weighted nearest-neighbour classifier. Conversely, Prototypical Networks \cite{snell2017prototypical} pose that an embedding exists around a single prototype representation for each class. It produces a distribution over classes for a query point based on a softmax of the distances to the prototypes. Both approaches have been successfully applied to vision and language tasks. 

In summary, we noticed that a GNN trained on a large dataset is the dominant tool used to learn to improve SOA MIP solvers. However, they require extensive training and tuning to yield acceptable performance.  In contrast, our method aims to be as efficient as possible by exploiting readily available data from the B\&B tree. As a consequence, we can take advantage of simpler models that are less prone to overfitting.  This design choice creates the opportunity for an accurate model with low opportunity costs.

\section{Methodology}

A CO problem $\P$ can be formulated as a mixed-integer linear programming (MIP) model using the notation
\begin{align}
     \P := \argmin_{\x} \{ c^T \x\,  |\,  A\x \leq b, \, \x\in \mathcal{B} \cup \mathcal{Q} \cup \mathcal{W}  \} ,  \label{eq:problem}
\end{align}
where $\mathbf{x} $ is the vector of decision variables that is partitioned into  $\mathcal{B}$, $\mathcal{Q}$ and $\mathcal{W}$, the sets of binary, integer and continuous variables, respectively. Our approach makes the assumption that some subset of the constraints are S1 \cite{beale1970special}.  The motivation for working with S1 constraints comes from the fact that it can be modelled as a classification problem. Indeed, S1 constraints impose that a sum of binary variables must be equal to 1, which implies that one variable must be chosen per constraint. This is equivalent to predicting the class in a classification problem. To link these concepts, we define a set of possible classes $K_v$ for constraint $v \in \V$, where each class $k \in K_v$ is associated with a different binary variable. This formulation allows us to use ML tools to predict the optimal class $\ko$ associated with each S1 constraint in problem $\P$. The solution vector for the binary variables is analogous to the one-hot encoding of the optimal class. The mathematical representation of S1 is given by 
\begin{align}
     \sum_{k \in K_v} x_v^k = 1   \hspace{1.5cm}   \forall v \in \V,  \label{eq:sos1}
\end{align}
where $x_v^k$ is the binary variable associated with the class $k \in K_v$. The set $\V$ contains all the S1 constraints where only a subset of them have their binary variables frozen. At inference, we generate a predicted class $\kt$ and we use it to freeze the corresponding binary variable $x_v^{k'}$.  We refer an S1 constraint $v$ as \emph{frozen} whenever we add the constraint $x_v^{\kt} = 1 $ in the MIP to create a reduced problem $\Pr$. \startblue Given that $\Pr$ is a smaller problem, we can expect it to be easier to solve than the original problem $\P$. Additionally, if a classifier can accurately predict the optimal class ($\kt = \ko $) for each constraint, then  $\Pr$ will contain the optimal solution. \stopblue The approach can be decomposed into three steps which are probe, select and freeze. These routines are explained in more detail in the remainder of this section. 

\subsection{Probe}
The probing step is used to collect data about the decision variables $\xvk$. We call the solver on $\P$ for a predefined probing time budget $\tp$ and fetch intermediate solutions in the callback. The probing routine does not discriminate any intermediate solution; it stores all integer feasible and fractional solutions. We define $\xvkt$ as the value of the variable $\xvk$ at iteration $t$, which corresponds to a node in the B\&B tree. At every iteration, we compute the most likely class based on which variable has the highest value $\kvt = \argmax_{k \in K_v} \{ \xvkt \}$.  We then concatenate the $\kvt$ values into the vector $\kp = [k_{v1}, k_{v2}, ..., k_{vn}]$, where $n$ is the probing sample size. The vector $\kp$  corresponds to the probing data for the S1 constraint $v$, and it is used as input to the classifier. Finally, we define $\K$ as the set of all probing data vectors $\kp$. 

\subsection{Select}
Given the heuristic nature of freezing variables, we need to select the S1 constraints whose variables will be the object of a freezing in a way that minimizes potential assignment errors. The selection strategy uses a scoring system to sort the S1 constraints based on the entropy $H$ to infer uncertainty. It is defined as 
\begin{align}
    H(\kp) &= - \sum_{k \in \kp} P(k | \kp) \log P(k| \kp) ,
\end{align}
where $P(k| \kp) $ is the probability associated with the class $k$ given the probing data $\kp$.  The probability is computed using the corresponding frequency:
\begin{align}        
    P(k| \kp) &= \frac{\left | \{ z \in \kp \, | \, z = k \} \right |}{\left | \kp \right |}. \label{eq:freq} 
\end{align}
Next, we compute the score which is the negative of the entropy:
\begin{align}
    \score(v) = - H(\kp).  \label{eq:score}
\end{align}
 Finally, we  select the first $\r \cdot \left | \V \right | $ constraints from the sorted list to produce the set of selected constraints $\Vr$. The ratio of constraints to freeze $\r$ is a hyperparameter that can be tuned to modify the aggressiveness of the algorithm.  We also experiment with a threshold-based strategy selecting constraints that have an entropy value below a given threshold.  The idea behind this strategy is reasonably intuitive. On the one hand, when the entropy $H(\kp)$ is low, it means that the solutions found during probing are similar and the most frequent class in $\kp$ is likely to be optimal. On the other hand, when the entropy is high, the most frequent class changes often and therefore we are less confident about what the optimal class is. 

\subsection{Freeze}
Once probing is done, we restart the solver and create a reduced problem $\Pr$ by freezing the variables $\xvk$ in the selected constraints $\Vr$. The freezing routine builds \emph{freezing cuts} ($FC$) defined as follows:
\begin{align}
    & FC(v, \kt):= x_v^{\kt} = 1, 
\end{align}
where the predicted class $\kt$ is the most likely class in the probing vector:
\begin{align}
    \kt = \argmax_k (P(k |  \kp)). 
\end{align}
The underlying  classifier uses the histogram method with discrete bins for each class in $\kp$. The reduced problem $\Pr$ is then solved to completion or until the B\&B process reaches the total time limit. The probe and the freeze methodology is summarized in Algorithm \ref{alg:pnf}. It is worth noting that both the select and the freeze routines are deterministic given some probing data. \startblue Therefore, the heuristic is not sensitive to the choice of the random seed. Futhermore, we demonstrate in Section~\ref{sec:results} that PNF is resilient to  different B\&B paths since it performs well even if probing is done using different solvers. \stopblue

\begin{algorithm}
    \caption{Probe and Freeze (PNF) algorithm}\label{alg:pnf}
    \begin{algorithmic}
    \Procedure{PNF}{$\P, \tp,\r$}
    \State $\Pr \gets \P$ \Comment{Make a copy of the problem}
    \State $\K \gets \probe ( P, \tp)$ \Comment{Collect probing data}
    \State $\Vr \gets \select (\K, \r)$ \Comment{Select the constraints to freeze}
    \For{$v \in \Vr$}
        \State $\kp \gets \K$ \Comment{Get probing data for the constraint $v$}
        \State $\kt \gets \argmax_k (P(k |  \kp))$ \Comment{Predict the optimal class}
        \State Add $FC(v, \kt)$ to $\Pr$  \Comment{Freeze the selected constraint}
    \EndFor
    \State $\x' \gets \solve(\Pr)$ \Comment{Solve the reduced problem and produce an approximate solution $\x'$}
    \EndProcedure
    \end{algorithmic}
\end{algorithm}

\textbf{Feasiblity considerations.} Since PNF fixes variables in different constraints independently, it is possible that the reduced problem is infeasible. In particular, the infeasibility issue can occur when the same variable appears in multiple S1 constraints and the predicted classes are not compatible. In practice, we can circumvent this issue by modifying the freezing cut to include both the predicted class and the class $k_{\text{inc}}$ from the incumbent solution:
\begin{align}
    & FC(v, \kt, k_{\text{inc}}):= x_v^{\kt} +  x_v^{k_{\text{inc}}} = 1. 
\end{align}
This modification ensures that the reduced problem is at least feasible with respect to the incumbent solution. \startblue 
This also circumvents the issue of conflicting assignments.\stopblue
If no incumbent solution is available from the probing phase,\startblue  PNF is still applicable, \stopblue but it cannot guarantee feasibility of the reduced problem. In case of infeasibility one can restart the solver with a different set of constraints. For example, the heuristic from \cite{hendel2022adaptive} dynamically changes the fixing rate based on the status of the current round. If no new solutions are found, the fixing rate is increased, whereas it is decreased when the reduced problem is infeasible or solved to optimality. \startblue For iterative algorithms, one can use an objective cut-off to exclude the incumbent and get out of a local optimum. However, this does not guarantee that the reduced problem is feasible. \stopblue \startblue 
In practice, we found that the infeasibility issue is not a major concern. As we report in the following section, for the LAP, the reduced problem is always feasible and, on MIPLIB instances, the infeasibility rate is less than 6.7\% (2/30). This is linked to the difficulty of finding feasible solutions for some MIPLIB instances. \stopblue

\section{Experimental Results} \label{sec:results}

We test our approach on two different sets of instances. The first set includes LAP instances constructed using historical data from CN. 
The second set contains instances with S1 constraints from the MIPLIB 2017 \cite{miplib2017} library. 
These two sets differ with respect to the underlying distribution. The first set is regular, whereas the second set is much more diverse. 
The LAP instances are generated using consecutive weeks of data from CN, which implies that all instances come from the same distribution.
In contrast, instances in MIPLIB were built independently and relate to many different problems. 
For this reason, MIPLIB can be seen as a more challenging test bed given our methodology targets a regular setting. 
Since each instance in MIPLIB is unique, it is not possible to effectively calibrate the two hyperparameters of PNF ($\r$, $\tp$).  
Nevertheless, we believe it is important to test PNF on MIPLIB to assess its limitations and understand where it may not be applicable. 
The results are presented in Sections~\ref{sec:lap} and \ref{sec:miplib}.

\textbf{Metrics. } The goal of this work is to quickly find high-quality solutions to MIPs with S1 constraints. Given that we are working with heuristics, the optimality gap as a performance metric is not appropriate because the optimal solution is not always known. Instead, we use the primal gap (PG) in percentage which is defined as follows:

\begin{align} \label{eq:RG}
    \text{PG} := \frac{\ct \x'- \ct \xc}{\ct \xc } \times 100, 
\end{align}
where $\x'$ is the best solution found by the heuristic and $\xc$  is the best known solution found across all experiments. We report several metrics: quantiles, average, and geometric average. The quantiles are useful to assess spread and skewness. The geometric mean is the $n^{\text{th}}$ root of the product of $n$ values, providing a measure that minimizes the effect of very large values. \startblue We add a shift factor of 1 in the calculation of the geometric mean to avoid issues with zero values. \stopblue Additionally, we count instances where the heuristic fails to find a feasible solution (No Sol) and instances where it outperforms every other algorithm (Wins).     \startblue To ensure a fair comparison, the runtime is calculated by adding the probing time budget to the solution process. We start the clock when the solution process starts, and we stop it when the best feasible solution is found.  Both RINS and LB use the best feasible solution found during probing. If no feasible solution is found, we report it as a failure (No Sol) because the heuristic cannot be executed. All heuristics have access to the same probing data and computational budget to maintain fairness. Given that the solver does not require probing, the probing time budget is added to its total time limit. \stopblue

\textbf{Setup. } Our software stack is fully written in Julia with the JuMP package \cite{lubin_jump_2023}. We compare the performance with the latest versions of three different solvers:  CPLEX version 22.11, Gurobi 10 and SCIP v8. By default, we run each scenario on a single thread of an Intel Gold 6148 Skylake (2.4 GHz) processor and 35 GiB of memory. 

\textbf{Baselines. } We compare the performance of PNF with a variety of heuristics.  We focus on heuristics that share the behavior of PNF, which is to create a reduced problem that should be easier to solve. The list of heuristics is as follows: 
\begin{itemize}
    \item \textbf{PNF: } The PNF heuristic fixes a percentage of the S1 constraints giving priority to the ones with the lowest entropy. 
    \item \textbf{PNFT: } Same as PNF but with a threshold instead of a fixing ratio. All constraints with an entropy value below the threshold are fixed. 
    \item \textbf{RINS: } The RINS heuristic \cite{danna2005exploring} fixes a percentage of the variables based on the matching set between the LP relaxation and the best feasible solution $\xinc$ from the probing phase. 
    \item \textbf{LB: } The LB heuristic \cite{fischetti_local_2003} creates a neighbourhood around the best feasible solution  $\xinc$ from the probing phase using a distance parameter $\dmax$. More formally, we add a constraint to the MIP to limit the binary norm relative to $\xinc$:  $\left\| {\x - \xinc} \right\|_b \leq \dmax  $. The value of $\dmax$ is calculated using a percentage of the number of binary variables in the MIP. For example, LB-0.8 would allow a maximum of 80\% of the variables to change relative to $\xinc$, which corresponds to a soft fixing ratio of 20\%.
    \item \textbf{GH: } Our Greedy Heuristic (GH) iteratively fixes the S1 constraints based on LP values $\xlp$. At every iteration $i$, it solves the LP relaxation and fixes the S1 constraint $v_i $ that is the least fractional: $v_i = \argmax_v \{\max_k \{\xvklp\}\}$. The process is repeated until we reach the given fixing ratio.  The main reason we include GH is to compare the performance of PNF with a heuristic that does not require a feasible solution to start.
 \end{itemize}
We implemented the heuristics RINS and LB following the descriptions in \cite{hendel2022adaptive}. Given that they are the best performing heuristics in the cited work that uses MIPLIB as a test set, there is evidence that they are high-quality baselines. \startblue Note that we also considered Relaxation Enforced Neighborhood Search (RENS) \cite{berthold_rens_2014} as a baseline, but we found that it is equivalent to RINS for binary variables. \stopblue All heuristics follow a similar naming convention where the prefix is the abbreviation and its hyperparameters are suffixes. Generally, the first suffix is a reference to the parameter that controls the size of the subproblem (fixing ratio, threshold or relative distance). The second suffix is the probing time budget $\tp$ in minutes. 

 \textbf{Tuning. } We tune the hyperparameters of the heuristics to produce conservative problem reductions. We opt for a standard setting across all scenarios to allow for a fair comparison. In our experiments, we test the 20\% to 50\% reduction range as it provides a good trade-off between performance and feasibility. Below 10\%, the reduction has a marginal effect, whereas a reduction above 50\% is likely to be too aggressive. We also use a standard setting for the probing time budget to ensure the same data is provided to each heuristic. The probing time budget is adjusted based on the difficulty of the instances in a given dataset. Some preliminary experiments were done to find reasonable values for the hyperparameters of each heuristic.

\subsection{LAP Instances}
\label{sec:lap}
The LAP is a scheduling problem that is used to assign locomotives to trains in a railway network. 
The problem is modelled using a space-time network in which each node represents a time and space location, whereas each arc represents a movement of locomotives.
The goal of the LAP is to find the cheapest way to send a certain number of locomotives through a network such that all capacity and demand constraints are satisfied. 
The main reasons that motivate the use of a learning heuristic for the LAP are the following.
First, the problem is NP-hard and it is challenging to find good solutions in a reasonable amount of time. 
Second, given the scheduling nature of the problem, instances follow daily and weekly cyclical patterns which are useful for learning. 

\textbf{MIP model. } The S1 constraints in the LAP are used to model the fact that only one set of locomotives (i.e., a consist) can be assigned to each train.
The decisions $ x_v^k$ represent the configuration of locomotives $k$ assigned to each train arc $v$ in the instance. 
The MIP of the LAP also contains flow conservation constraints for each node in the network. We will not discuss the details of the model here as it is not the focus of this paper.
We refer to \cite{camilo} for the details on the MIP formulation of the LAP. 

\textbf{Datasets: LAP-A and LAP-B.} The historical data used to build the instances were provided by CN. The dataset contains the origin and destination stations for each train operated during the year. 
We study the performance on different difficulty levels to access the generalization properties of the heuristic. We run our experiments on two datasets with different difficulty levels: LAP-A and LAP-B. The LAP-A dataset is composed of 30 instances with 336 trains per instance on average. The instances in LAP-A are generated by a subsample process that preserves the weekly and daily patterns. 
This is done by selecting the set of trains to operate in a given region of the network. We tune the difficulty level such that instances can be solved within a reasonably short amount of time (less than 1 hour). The goal of LAP-A is to test many different configurations of heuristics in an efficient manner. The second dataset LAP-B is composed of 10 instances with 5673 trains on average. They contain all trains that operate in the network in a given week. The goal of LAP-B is to test the performance of the heuristic on the most difficult setting. The LAP-B instances are generated using the same methodology as in \cite{camilo} and we run them with a time limit of 6 hours. The main difference compared to \cite{camilo} is that we test on the most difficult instances (all trains) whereas \cite{camilo} split the test set into mainline only and all trains. 

\textbf{Discussion on probing data. } Before discussing the performance metrics, we provide some insights on the probing data. The probing data is key to our method since we rely on it to select which variables to freeze. In particular, we are interested in the number of samples available, which is directly linked to the number of B\&B nodes visited during the probing phase. We show the number of nodes visited for LAP-A and LAP-B in Tables \ref{p_lap_a} and \ref{p_lap_b}, respectively. For LAP-A, CPLEX visits 1,382 nodes on average, whereas SCIP only visits 20. Gurobi visits around 900 nodes in the same amount of time. We observe a similar pattern for LAP-B. The number of nodes visited by CPLEX is almost 80 times higher than SCIP on LAP-B. This is a significant difference that can explain the performance gap between the two solvers.

\textbf{Results. } The primal gap for the LAP-A dataset is presented for each solver in Tables \ref{tab:gap_cplex_lap_a}, \ref{tab:gap_gurobi_lap_a} and \ref{tab:gap_scip_lap_a}. For this dataset, the probing time budget is 2 minutes and the time limit is 10 minutes overall. For CPLEX, PNF-0.5-2 dominates with the lowest average gap of 0.23\% and 12 wins out of 30. The second best scenario is PNFT-0.05-2 with 9 wins followed by CPLEX alone with 7 wins.  For Gurobi, the best scenario is PNFT-0.0-2 with an average gap of 0.03\% and 13 wins. Both PNF configurations perform well with 12 wins each. The feasibility focus mode (grb-f1) has a marginal effect on LAP-A because it is relatively easy to find feasible solutions for these instances. We note that LB is the next best heuristic for Gurobi but it is not competitive with CPLEX. For SCIP, the scenario PNFT-0.05-2 is also the best with an average gap of 0.96\% and 20 wins. We observe a large number of No Sol for LB and RINS because SCIP fails to find a feasible solution for half the instances during the probing phase. A feasible solution is essential for the execution of both heuristics. We found that GH is not competitive for any solver. 

In Table~\ref{tab:time_cplex_lap_a}, we show the statistics for the runtime using CPLEX. The RINS-0.5-2 scenario is significantly ahead because it is the first to find its best solution for 26 out of 30 instances. However, it has the worst average gap at 1.39\%. In this setting, the PNF heuristic does not have a meaningful edge in terms of runtime. This is because the time limit is relatively short at 10 minutes such that most runs terminate early. We found similar runtime values for all solvers.

In Table~\ref{tab:fixing_ratio_cplex_lap_a}, we show the equivalent fixing ratio for PNFT where the constraint selection is done using a threshold. The threshold at 0 implies that we only fix constraints where all the solutions from the probing data agree on the value of the variables. On average, this threshold is equivalent to fixing 26\% of the constraints. The threshold at 0.05 is equivalent to fixing 36\% of the constraints. Also, we observe that the fixing ratio can vary significantly across instances (from 19\% to 50\%). This is because the value of entropy for each constraint varies across instances. 

The results for the LAP-B dataset are presented to highlight the ability of PNF to find good feasible solutions for large-scale instances. We report fewer scenarios for this dataset given its large computational burden. The time limit is 6 hours for LAP-B. The primal gaps for CPLEX and Gurobi are reported in Tables~\ref{tab:gap_cplex_lap_b} and \ref{tab:gap_gurobi_lap_b}, respectively. The LAP-B instances benefit from a higher fixing ratio since PNF-0.5-30 is the best scenario for both solvers. The LB heuristic is competitive on smaller LAP instance, but it is not able to find any solution with CPLEX. The results on Gurobi are similar where LB is not able to find any solution for 3 out of 10 instances. 

For LAP-B, we also compute the performance over time as shown in Tables~\ref{tab:gap_over_time} and \ref{tab:sol_found_over_time}. These tables respectively display the gap relative to the best known solution and the number of instances solved at different points in time. On average, with the best configuration (PNF-0.5-30), PNF is able to find an upper bound within 1.48\% of the best known solution after 1 hour (30 minutes for probing + 30 minutes for solving) whereas CPLEX's average gap is above 4\%. Furthermore, PNF finds a feasible upper bound for more instances during the first two hours (10/20 for PNF-0.2-30 and 7/20 for CPLEX). 

\textbf{Key findings.} The results show that our method is the most effective one on the LAP. The version with a fixing ratio (PNF) performs best on CPLEX, whereas the threshold-based version (PNFT) is preferred on the two other solvers. RINS is often the fastest but produces low-quality solutions. LB is competitive on LAP-A but it can fail to find any solution on LAP-B. GH is not competitive for any solver.


\begin{table}[h]
\centering
\caption{B\&B node count from probing 2 minutes on LAP-A}
\label{p_lap_a}
\begin{tabular}{lcccccc}
\toprule
{Scenario} & \multicolumn{3}{c}{Quantiles} & {Mean} & {Geomean} \\
{} & {0.1} & {0.5} & {0.9} & {} & {} \\
\midrule
grb-probing & 137.70 & 548.50 & 2105.30 & 899.87 & 536.73 \\
cplex-probing & 612.60 & 1142.50 & 3114.40 & 1381.93 & 1152.45 \\
scip-probing & 17.00 & 20.00 & 25.00 & 20.33 & 21.09 \\
\bottomrule
\end{tabular}
\end{table}

\begin{table}[h]
\centering
\caption{Primal gap (\%) using CPLEX with a 10-minute limit on LAP-A}
\label{tab:gap_cplex_lap_a}
\begin{tabular}{lcccccccc}
\toprule
{Scenario} & \multicolumn{3}{c}{Quantiles} & {Mean} & {Geomean} & {No Sol} & {Wins} \\
{} & {0.1} & {0.5} & {0.9} & {} & {} & {} & {} \\
\midrule
PNFT-0.05-2 & 0.00 & 0.06 & 0.64 & 0.28 & 1.23 & 0 & 9 \\
PNFT-0.0-2 & 0.00 & 0.18 & 1.17 & 0.39 & 1.32 & 0 & 4 \\
PNF-0.2-2 & 0.00 & 0.11 & 1.27 & 0.43 & 1.35 & 0 & 6 \\
PNF-0.5-2 & 0.00 & 0.06 & 0.81 & 0.23 & 1.20 & 0 & 12 \\
RINS-0.2-2 & 0.10 & 0.57 & 1.90 & 0.78 & 1.67 & 1 & 0 \\
RINS-0.5-2 & 0.19 & 1.14 & 3.04 & 1.39 & 2.16 & 1 & 0 \\
LB-0.7-2 & 0.00 & 0.35 & 1.38 & 0.58 & 1.47 & 1 & 1 \\
LB-0.8-2 & 0.00 & 0.28 & 1.87 & 0.62 & 1.48 & 1 & 1 \\
cplex & 0.00 & 0.34 & 1.96 & 0.70 & 1.54 & 0 & 7 \\
GH-0.05 & 0.00 & 0.20 & 1.93 & 0.65 & 1.50 & 0 & 3 \\
GH-0.1 & 0.00 & 0.24 & 1.44 & 0.57 & 1.45 & 0 & 2 \\
\bottomrule
\end{tabular}
\end{table}

\begin{table}[h]
\centering
\caption{Primal gap (\%) using Gurobi with a 10-minute limit on LAP-A}
\label{tab:gap_gurobi_lap_a}
\begin{tabular}{lcccccccc}
\toprule
{Scenario} & \multicolumn{3}{c}{Quantiles} & {Mean} & {Geomean} & {No Sol} & {Wins} \\
{} & {0.1} & {0.5} & {0.9} & {} & {} & {} & {} \\
\midrule
PNFT-0.05-2 & 0.00 & 0.00 & 0.07 & 0.03 & 1.03 & 0 & 10 \\
PNFT-0.0-2 & 0.00 & 0.00 & 0.07 & 0.03 & 1.03 & 0 & 13 \\
PNF-0.2-2 & 0.00 & 0.00 & 0.26 & 0.08 & 1.07 & 0 & 12 \\
PNF-0.5-2 & 0.00 & 0.00 & 0.12 & 0.03 & 1.03 & 0 & 12 \\
RINS-0.2-2 & 0.00 & 0.08 & 0.46 & 0.16 & 1.15 & 0 & 3 \\
RINS-0.5-2 & 0.00 & 0.09 & 0.64 & 0.26 & 1.21 & 0 & 3 \\
LB-0.7-2 & 0.00 & 0.01 & 0.21 & 0.07 & 1.06 & 0 & 10 \\
LB-0.8-2 & 0.00 & 0.01 & 0.21 & 0.08 & 1.07 & 0 & 10 \\
grb-f1 & 0.00 & 0.01 & 0.15 & 0.06 & 1.05 & 0 & 8 \\
grb & 0.00 & 0.01 & 0.30 & 0.11 & 1.09 & 0 & 9 \\
GH-0.05 & 0.00 & 0.07 & 0.50 & 0.16 & 1.14 & 0 & 6 \\
GH-0.1 & 0.00 & 0.08 & 0.43 & 0.20 & 1.16 & 0 & 5 \\
\bottomrule
\end{tabular}
\end{table}

\begin{table}[h]
\centering
\caption{Primal gap (\%) using SCIP with a 10-minute limit on LAP-A}
\label{tab:gap_scip_lap_a}
\begin{tabular}{lcccccccc}
\toprule
{Scenario} & \multicolumn{3}{c}{Quantiles} & {Mean} & {Geomean} & {No Sol} & {Wins} \\
{} & {0.1} & {0.5} & {0.9} & {} & {} & {} & {} \\
\midrule
PNFT-0.05-2 & 0.14 & 0.68 & 2.05 & 0.96 & 1.81 & 0 & 20 \\
PNFT-0.0-2 & 0.14 & 0.68 & 2.32 & 1.00 & 1.83 & 0 & 18 \\
PNF-0.2-2 & 0.81 & 2.98 & 5.03 & 3.13 & 3.74 & 0 & 1 \\
PNF-0.5-2 & 0.04 & 1.74 & 3.90 & 2.00 & 2.59 & 0 & 6 \\
RINS-0.2-2 & 0.58 & 0.99 & 2.51 & 1.24 & 2.12 & 15 & 0 \\
RINS-0.5-2 & 2.64 & 4.41 & 4.99 & 4.10 & 5.00 & 15 & 0 \\
LB-0.7-2 & 0.98 & 2.88 & 5.61 & 3.17 & 3.70 & 15 & 0 \\
LB-0.8-2 & 0.65 & 2.88 & 5.29 & 3.11 & 3.63 & 15 & 0 \\
scip & 1.59 & 3.27 & 4.94 & 3.35 & 4.05 & 0 & 1 \\
GH-0.05 & 1.34 & 3.05 & 4.26 & 3.03 & 3.75 & 0 & 0 \\
GH-0.1 & 0.79 & 2.96 & 4.53 & 2.87 & 3.53 & 0 & 2 \\
\bottomrule
\end{tabular}
\end{table}

\begin{table}[h]
\centering
\caption{Runtime (s) using CPLEX with a 10-minute limit on LAP-A}
\label{tab:time_cplex_lap_a}
\begin{tabular}{lccccccc}
\toprule
{Scenario} & \multicolumn{3}{c}{Quantiles} & {Mean} & {Geomean} & {Wins} \\
{} & {0.1} & {0.5} & {0.9} & {} & {} & {} \\
\midrule
PNFT-0.05-2 & 198.12 & 536.75 & 596.76 & 477.51 & 446.04 & 2 \\
PNFT-0.0-2 & 217.93 & 564.10 & 601.13 & 480.78 & 450.85 & 0 \\
PNF-0.2-2 & 258.25 & 480.52 & 600.19 & 454.41 & 426.45 & 0 \\
PNF-0.5-2 & 160.38 & 472.92 & 598.48 & 416.94 & 371.21 & 0 \\
RINS-0.2-2 & 122.06 & 241.79 & 588.52 & 343.14 & 279.69 & 2 \\
RINS-0.5-2 & 121.22 & 121.32 & 123.11 & 125.79 & 125.63 & 26 \\
LB-0.7-2 & 332.50 & 553.67 & 592.00 & 479.91 & 456.60 & 0 \\
LB-0.8-2 & 238.17 & 496.44 & 592.29 & 454.67 & 424.94 & 0 \\
cplex & 209.44 & 510.97 & 589.85 & 420.73 & 378.57 & 0 \\
GH-0.05 & 204.64 & 532.53 & 593.81 & 463.06 & 428.56 & 0 \\
GH-0.1 & 208.88 & 553.53 & 590.87 & 449.79 & 400.71 & 2 \\
\bottomrule
\end{tabular}
\end{table}

\begin{table}[h]
\centering
\caption{Fixing ratio with CPLEX on LAP-A}
\label{tab:fixing_ratio_cplex_lap_a}
\begin{tabular}{lcccccc}
\toprule
{Scenario} & \multicolumn{3}{c}{Quantiles} & {Mean} & {Geomean} \\
{} & {0.1} & {0.5} & {0.9} & {} & {} \\
\midrule
PNFT-0.05-2 & 0.19 & 0.39 & 0.50 & 0.36 & 1.36 \\
PNFT-0.0-2 & 0.12 & 0.28 & 0.39 & 0.26 & 1.25 \\
\bottomrule
\end{tabular}
\end{table}

\begin{table}[h]
\centering
\caption{B\&B node count from probing 30 minutes on LAP-B}
\label{p_lap_b}
\begin{tabular}{lcccccc}
\toprule
{Scenario} & \multicolumn{3}{c}{Quantiles} & {Mean} & {Geomean} \\
{} & {0.1} & {0.5} & {0.9} & {} & {} \\
\midrule
grb-probing & 168.10 & 362.00 & 384.40 & 295.00 & 269.17 \\
scip-probing & 7.60 & 13.50 & 28.00 & 16.70 & 15.33 \\
cplex-probing & 1038.00 & 1263.00 & 1610.50 & 1311.00 & 1284.92 \\
\bottomrule
\end{tabular}
\end{table}

\begin{table}[h]
\centering
\caption{Primal gap (\%) using CPLEX with a 6-hour limit on LAP-B}
\label{tab:gap_cplex_lap_b}
\begin{tabular}{lcccccccc}
\toprule
{Scenario} & \multicolumn{3}{c}{Quantiles} & {Mean} & {Geomean} & {No Sol} & {Wins} \\
{} & {0.1} & {0.5} & {0.9} & {} & {} & {} & {} \\
\midrule
PNFT-0.1-30 & 0.12 & 0.26 & 1.08 & 0.60 & 1.49 & 0 & 2 \\
PNF-0.2-30 & 0.52 & 1.05 & 1.63 & 1.25 & 2.08 & 0 & 1 \\
PNF-0.5-30 & 0.00 & 0.15 & 1.22 & 0.35 & 1.28 & 0 & 7 \\
LB-0.7-30 & nan & nan & nan & nan & nan & 10 & 0 \\
cplex & 0.58 & 2.68 & 5.40 & 4.51 & 3.73 & 0 & 0 \\
\bottomrule
\end{tabular}
\end{table}

\begin{table}[h]
\centering
\caption{Primal gap (\%) using Gurobi with a 6-hour limit on LAP-B}
\label{tab:gap_gurobi_lap_b}
\begin{tabular}{lcccccccc}
\toprule
{Scenario} & \multicolumn{3}{c}{Quantiles} & {Mean} & {Geomean} & {No Sol} & {Wins} \\
{} & {0.1} & {0.5} & {0.9} & {} & {} & {} & {} \\
\midrule
PNFT-0.1-30 & 0.00 & 0.15 & 0.33 & 0.18 & 1.17 & 0 & 3 \\
PNF-0.2-30 & 0.00 & 4.73 & 9.02 & 4.59 & 3.77 & 0 & 1 \\
PNF-0.5-30 & 0.00 & 0.07 & 0.27 & 0.11 & 1.10 & 0 & 5 \\
LB-0.7-30 & 0.11 & 10.22 & 10.98 & 7.45 & 5.82 & 3 & 1 \\
grb & 4.05 & 5.03 & 5.43 & 4.86 & 5.82 & 0 & 0 \\
grb-f1 & 0.50 & 1.05 & 4.70 & 2.19 & 2.42 & 0 & 0 \\
\bottomrule
\end{tabular}
\end{table}

\begin{table}[h]
\footnotesize 
\begin{minipage}[b]{0.5\textwidth}
\centering
\caption{Primal gap over time for LAP-B}
\label{tab:gap_over_time}
\begin{tabular}{lccccc}
\toprule
{Scenario} & \multicolumn{4}{c}{Time after probing (s)} \\
{} & {1800} & {3600} & {5400} & {7200} \\
\midrule
cplex & nan & 4.04 & 3.94 & 2.49 \\
PNF-0.2-10 & 1.72 & 1.48 & 1.09 & 0.88 \\
PNF-0.2-30 & 2.01 & 1.67 & 1.18 & 1.17 \\
PNF-0.5-30 & 1.48 & 1.39 & 0.89 & 0.69 \\
\bottomrule
\end{tabular}

\end{minipage}
\begin{minipage}[b]{0.5\textwidth}
\centering
\caption{LAP-B instances with feasible solution }
\label{tab:sol_found_over_time}
\begin{tabular}{lccccc}
\toprule
{Scenario} & \multicolumn{4}{c}{Time after probing (s)} \\
{} & {1800} & {3600} & {5400} & {7200} \\
\midrule
cplex & 0 & 6 & 7 & 7 \\
PNF-0.2-10 & 7 & 9 & 9 & 9 \\
PNF-0.2-30 & 7 & 10 & 10 & 10 \\
PNF-0.5-30 & 8 & 9 & 9 & 9 \\
\bottomrule
\end{tabular}

\end{minipage} 
\end{table}

\subsection{MIPLIB Instances}
\label{sec:miplib}
The MIPLIB \cite{miplib2017} dataset is an open source library of MIP instances. It contains a variety of instances from different application domains such as transportation, scheduling, and others related to combinatorial optimization. 
The instances in MIPLIB are generated independently and are not related to each other.
This makes it a more challenging test bed for our approach.
The main reason for using MIPLIB is to show that our approach can generalize to instances that have S1 constraints as the only common trait. \startblue This is challenging because there is a meaningful amount of complexity in each instance that our heuristic does not capture since it only considers variables that appear within S1 constraints. Furthermore, the number and the types of constraints can vary significantly from one instance to the other. \stopblue
In that sense, the PNF heuristic is unlikely to perform well on all MIPLIB instances. 
The task of selecting the appropriate algorithm for a given problem is often called the Algorithm Selection Problem (ASP)~\cite{algorithmselection}. 
This can be done by building a mapping between the features of the problem and the expected performance of an algorithm. 
In our case, we are interested in selecting the MIPLIB instances on which the PNF heuristic is likely to perform well. 
To achieve this, we identify the features related to the implicit assumption of the approach. 

\textbf{Dataset: MIPLIB-S1.}  The main assumption is related to the presence of S1 constraints in the MIP. From the initial set of more than a thousand instances, we select the ones that contain at least one S1 constraint. After this first filtering step, we are left with 188 instances. The second filtering step aims to select instances that strike a balance in complexity.  On the one hand, if the instance can be solved to optimality or close to optimality in a short amount of time, then the PNF heuristic is not needed. On the other hand, if the instance is too difficult, then there will not be enough probing data to make good decisions. It is well known that the quantity of data is key to reach high performance for ML techniques. To achieve this balance, we do a probing run on each instance to measure both the optimality gap and the number of nodes visited. Instances with fewer than 10 visited nodes are removed because the sample size is too small to make reliable predictions. We then select 30 instances with the highest optimality gap and remove outliers with a gap above a thousand. This is because primal heuristics are more likely to be effective on instances where the optimality gap is difficult to close.\startblue We decided to use a dataset size of 30 because of the limited number of instances available after the filtering process. Also, it is considered sufficient to be statistically relevant.  \stopblue 

We run the selection process using different settings to build two datasets, MIPLIB-S1-A and MIPLIB-S1-B, to evaluate our performance for different difficulty levels. The MIPLIB-S1-A contains instances of medium difficulty and MIPLIB-S1-B contains instances of hard difficulty. The first dataset is built using a probing time budget of 2 minutes with the SCIP solver. The second dataset is built using a probing time budget of 1 hour with the Gurobi solver. Given that preliminary findings reveal that Gurobi is the top-performing solver, we classify instances that pose a challenge to it as difficult. However, SCIP typically makes little to no meaningful progress on many of them. The construction of MIPLIB-S1-A with SCIP ensures that every instance in that dataset is tractable for every solver.  We should mention that the notion of difficulty is not absolute but relative to a given solver, as some instances are more difficult for a given solver than another. Therefore, it is not obvious how to split MIPLIB instances into different difficulty levels in a way that is consistent across all solvers. In the end, we believe that the two datasets we propose offer a good compromise in that regard. 


    
\textbf{Results.} 
We follow a similar structure as in Section~\ref{sec:lap} to present the results for MIPLIB. We split the results for the same three solvers in Tables~\ref{tab:gap_cplex_mip_a}, \ref{tab:gap_gurobi_mip_a} and \ref{tab:gap_scip_mip_a}. The probing time budget is set to 2 minutes and 5 minutes for MIPLIB-S1-A and MIPLIB-S1-B, respectively. We allow a total time limit of 10 minutes for MIPLIB-S1-A and 20 minutes for MIPLIB-S1-B.  The greedy heuristic was removed from the comparison because of technical issues: solving the LP repeatedly is either too slow or leads to numerical issues for many instances. 

In this discussion, we refer to the results on MIPLIB-S1-A but the results on MIPLIB-S1-B lead us to the same conclusions. The table of results for MIPLIB-S1-B is included in Appendix~\ref{sec:appendix}. Overall, PNF excels when paired with CPLEX but it is relatively less successful when used with either of the other two solvers.  The specific reason for the performance disparity among the solvers cannot be determined with certainty because each one has its unique implementation of B\&B, which is not always visible to us. We suspect that the reason for CPLEX's better performance is its ability to create significantly more probing data than other solvers. This is consistent with our results on the LAP. For CPLEX, the best scenario is PNFT-0.0-2 with an average gap of 7.17\% and 19 wins. The second best scenario is PNFT-0.2-2 with 17 wins followed by CPLEX alone with 16 wins. With Gurobi, we found that LB-0.7 performs best with an average gap close to 2\% and 18 wins. The PNF heuristic outperforms RINS with this solver. Lastly, SCIP outperforms all the heuristics we tested because it has the highest number of wins. The scenario RINS-0.2-2 has the lowest average gap of 12.92\% but it does not find a feasible solution for three out of 30 instances. We qualify the performance of PNF with SCIP as neutral because it is similar to the different heuristics. This can be explained by the relatively low number of B\&B nodes generated by SCIP. As a final point, we highlight the fact that properly tuned PNF and PNFT are rarely the worst performing heuristic. The 90$^{\text{th}}$ quantile of the gap distribution is not higher than the alternatives. This is a good indication that they are safe to use compared to the other heuristics we tested.

\textbf{Key insights.} The choice of heuristic and solver combination is crucial in MIP problem solving. PNF works well with CPLEX, while Gurobi and SCIP have other preferable options. We found that no heuristic systematically outperformed the others. We believe this is a consequence of the heterogeneous nature of the MIPLIB dataset. The instances are not related to each other and, therefore, it is difficult to find a heuristic that performs well on all of them. This is in contrast with the LAP where the instances are similar to each other and, therefore, we would expect a heuristic that performs well on one instance to perform well on the others.

\begin{table}[h]
\centering
\caption{Primal gap (\%) using CPLEX with a 10-minute limit on MIPLIB-S1-A}
\label{tab:gap_cplex_mip_a}
\begin{tabular}{lcccccccc}
\toprule
{Scenario} & \multicolumn{3}{c}{Quantiles} & {Mean} & {Geomean} & {No Sol} & {Wins} \\
{} & {0.1} & {0.5} & {0.9} & {} & {} & {} & {} \\
\midrule
PNFT-0.2-2 & 0.00 & 0.00 & 16.89 & 7.39 & 2.81 & 0 & 17 \\
PNFT-0.0-2 & 0.00 & 0.00 & 16.20 & 7.17 & 2.67 & 0 & 19 \\
PNF-0.2-2 & 0.00 & 0.13 & 38.05 & 18.57 & 3.61 & 0 & 14 \\
PNF-0.5-2 & 0.00 & 0.12 & 29.78 & 15.06 & 3.56 & 0 & 14 \\
RINS-0.2-2 & 0.00 & 0.31 & 34.16 & 12.84 & 3.49 & 3 & 13 \\
RINS-0.5-2 & 0.00 & 0.24 & 36.73 & 18.34 & 4.24 & 3 & 13 \\
LB-0.7-2 & 0.00 & 0.00 & 19.82 & 7.34 & 2.33 & 3 & 15 \\
LB-0.8-2 & 0.00 & 0.00 & 23.33 & 9.32 & 2.90 & 3 & 14 \\
cplex & 0.00 & 0.00 & 24.48 & 8.07 & 2.76 & 1 & 16 \\
\bottomrule
\end{tabular}
\end{table}

\begin{table}[h]
\centering
\caption{Primal gap (\%) using Gurobi with a 10-minute limit on MIPLIB-S1-A}
\label{tab:gap_gurobi_mip_a}
\begin{tabular}{lcccccccc}
\toprule
{Scenario} & \multicolumn{3}{c}{Quantiles} & {Mean} & {Geomean} & {No Sol} & {Wins} \\
{} & {0.1} & {0.5} & {0.9} & {} & {} & {} & {} \\
\midrule
PNFT-0.2-2 & 0.00 & 0.00 & 10.68 & 2.98 & 1.91 & 2 & 12 \\
PNFT-0.0-2 & 0.00 & 0.00 & 9.14 & 2.30 & 1.81 & 2 & 12 \\
PNF-0.2-2 & 0.00 & 0.00 & 7.96 & 2.33 & 1.73 & 2 & 14 \\
PNF-0.5-2 & 0.00 & 0.00 & 8.56 & 2.78 & 1.90 & 2 & 11 \\
RINS-0.2-2 & 0.00 & 0.00 & 10.68 & 2.89 & 1.97 & 2 & 10 \\
RINS-0.5-2 & 0.00 & 0.00 & 16.28 & 5.00 & 2.36 & 2 & 11 \\
LB-0.7-2 & 0.00 & 0.00 & 2.72 & 2.02 & 1.36 & 2 & 18 \\
LB-0.8-2 & 0.00 & 0.00 & 2.53 & 1.99 & 1.38 & 2 & 16 \\
grb-f1 & 0.00 & 0.00 & 8.88 & 14.27 & 1.82 & 0 & 15 \\
grb & 0.00 & 0.00 & 7.80 & 14.83 & 1.74 & 0 & 17 \\
\bottomrule
\end{tabular}
\end{table}

\begin{table}[h]
\centering
\caption{Primal gap (\%) using SCIP with a 10-minute limit on MIPLIB-S1-A}
\label{tab:gap_scip_mip_a}
\begin{tabular}{lcccccccc}
\toprule
{Scenario} & \multicolumn{3}{c}{Quantiles} & {Mean} & {Geomean} & {No Sol} & {Wins} \\
{} & {0.1} & {0.5} & {0.9} & {} & {} & {} & {} \\
\midrule
PNFT-0.2-2 & 0.00 & 0.91 & 61.85 & 17.74 & 4.99 & 0 & 11 \\
PNFT-0.0-2 & 0.00 & 0.83 & 60.86 & 17.87 & 4.80 & 0 & 13 \\
PNF-0.2-2 & 0.00 & 6.90 & 48.84 & 20.67 & 6.38 & 0 & 11 \\
PNF-0.5-2 & 0.00 & 5.64 & 38.14 & 26.64 & 7.09 & 0 & 11 \\
RINS-0.2-2 & 0.00 & 0.86 & 33.97 & 12.92 & 3.95 & 3 & 13 \\
RINS-0.5-2 & 0.00 & 5.40 & 53.54 & 25.09 & 6.36 & 3 & 9 \\
LB-0.7-2 & 0.00 & 0.62 & 72.46 & 30.55 & 4.75 & 0 & 14 \\
LB-0.8-2 & 0.00 & 3.07 & 74.40 & 32.77 & 6.37 & 0 & 13 \\
scip & 0.00 & 0.19 & 46.91 & 22.80 & 4.36 & 0 & 16 \\
\bottomrule
\end{tabular}
\end{table}

\section{Perspectives and Future Work}
\label{sec:conclusion}

The main purpose of the paper was to explore if a one-shot learning heuristic could improve the performance of a commercial solver on MIPs with S1 constraints.  
According to our reading of the literature, most other works on this subject use some form of neural networks combined with supervised or reinforcement learning. The overhead of DL is problematic in practice because of the specialized hardware requirements and the challenge of generating representative training datasets. Our method bypasses this issue by using a classical ML model that can be trained using readily available data. This means that it can be integrated natively within a solver without requiring any additional hardware or data. \startblue The requirements to integrate our method are minimal. The solver must keep track of intermediate solutions in the B\&B tree and compute the entropy for each S1 constraint. PNF can then be used as a drop-in replacement for RINS or other similar heuristics. \stopblue 

For both small and large LAP instances, we were able to outperform alternative heuristics with three different solvers. On the MIPLIB dataset, the results were mixed. Only CPLEX enabled us to consistently outperform the other heuristics. This is likely due to the fact that CPLEX is the solver that generates the largest amount of data during the probing phase. We also noticed that the LB heuristic was most effective with Gurobi which highlights the importance of selecting the right heuristic for the right solver. 

This study has demonstrated that the initial solutions found during the first few iterations of B\&B can serve as reliable proxies for the optimal solution.  A significant contribution of this research lies in the use of entropy \startblue to make good heuristic decisions in the context of B\&B methods. By identifying and avoiding freezing the constraints with high entropy, we aim to limit the risk of compromising the quality of the final solution. \stopblue However, we found that our heuristic is less effective when the amount of probing data is limited. This can occur when using the default B\&B in SCIP. We suggest researchers that use SCIP to explore different ways to generate probing data. For example, one could iteratively solve the LP relaxation with custom perturbations to explore the solution space more efficiently. 
 
Future research might attempt to augment the heuristic by adding more features. However, we found that this did not improve the performance enough to justify the additional cost. In fact, we realized that the simpler heuristic was more robust and performed better in practice. 
This is consistent with Occam's razor principle, which states that the simplest explanation is usually the best. For that reason, we believe PNF should be seen as a strong baseline for future works on learned heuristics for MIPs because it is easy to understand and reproduce. 
An interesting next step would be to add a worst-case analysis and a feasibility recovery strategy on top of the heuristic to increase its robustness and reliability. 


\section*{Statements and Declarations}
\newcommand{\statement}[1]{\subsection*{#1}}

\statement{Acknowledgments}
We acknowledge the support of the Institute for Data Valorization (IVADO), Digital Alliance Canada, the Canada First Research Excellence Fund (Apogée/CFREF), and the Canadian National Railway Company (through the CN Chain in Optimization of Railway Operations) for their financial backing and provision of valuable data which was  instrumental in advancing this research. We are grateful to two anonymous reviewers for their valuable comments, which helped improve the quality of the paper. 

\statement{Ethical approval}
No ethical approval was required for this study.

\statement{Competing interests}
We, the authors of this study, declare that we have no competing interests. This encompasses both financial and non-financial interests that could undermine the integrity of our research.

\statement{Authors' contributions}
Charly Robinson La Rocca was involved in all aspects of the research, conceptualization, methodology, programming, preparing experimental results, and writing the manuscript. Jean-Fran{\c{c}}ois Cordeau and Emma Frejinger contributed to the conceptualization, methodology, and editing the manuscript. All authors reviewed the manuscript.

\statement{Funding}
Charly Robinson La Rocca received funding from Institute for Data Valorization (IVADO). 
Emma Frejinger received funding from the Canada Research Chairs Program. 

\statement{Availability of data and materials}
The data that support the findings of this study come from the Canadian National Railway Company. The data for the LAP is subject to confidentiality restrictions. However, the code and the related data for MIPLIB instances are open-source, accessible at \url{https://github.com/laroccacharly/MLSOS}. 

\newpage 
\begin{appendices}

\section{MIPLIB-S1-B Results}\label{sec:appendix}
This appendix contains the results for the MIPLIB-S1-B dataset. 
We report the primal gap for CPLEX, Gurobi and SCIP in Tables~\ref{tab:gap_cplex_mip_b}, \ref{tab:gap_gurobi_mip_b} and \ref{tab:gap_scip_mip_b}, respectively.
The Tables~\ref{tab:gap_time_cplex_miplib_b} and \ref{tab:primal_gap_time_grb_miplib_b} shows for each instance the primal gap and runtime for CPLEX and Gurobi, respectively. \startblue Finally, Table~\ref{tab:primal_integral} shows the primal integral (PI) and the number of solutions found by the best algorithms on this dataset. The minimum value in each column is underlined and highlighted in bold font. \stopblue 

\end{appendices}

\begin{table}[h]
\centering
\caption{Primal gap (\%) using CPLEX with a 20-minute limit on MIPLIB-S1-B}
\label{tab:gap_cplex_mip_b}
\begin{tabular}{lcccccccc}
\toprule
{Scenario} & \multicolumn{3}{c}{Quantiles} & {Mean} & {Geomean} & {No Sol} & {Wins} \\
{} & {0.1} & {0.5} & {0.9} & {} & {} & {} & {} \\
\midrule
PNFT-0.2-5 & 0.00 & 3.76 & 66.32 & 24.27 & 6.83 & 2 & 14 \\
PNFT-0.0-5 & 0.00 & 7.37 & 66.32 & 23.67 & 7.41 & 2 & 12 \\
PNF-0.2-5 & 0.00 & 4.51 & 62.86 & 19.44 & 6.46 & 1 & 14 \\
PNF-0.5-5 & 0.00 & 7.32 & 63.16 & 23.76 & 7.72 & 3 & 8 \\
RINS-0.2-5 & 0.00 & 7.79 & 100.07 & 82.69 & 7.94 & 5 & 12 \\
RINS-0.5-5 & 0.00 & 9.53 & 114.48 & 44.02 & 10.41 & 5 & 9 \\
LB-0.7-5 & 0.00 & 8.65 & 115.38 & 177.36 & 10.48 & 4 & 9 \\
LB-0.8-5 & 0.00 & 5.15 & 99.27 & 116.42 & 8.00 & 4 & 8 \\
cplex & 0.00 & 8.46 & 91.31 & 568.14 & 9.94 & 2 & 9 \\
\bottomrule
\end{tabular}
\end{table}

\begin{table}[h]
\centering
\caption{Primal gap (\%) using Gurobi with a 20-minute limit on MIPLIB-S1-B}
\label{tab:gap_gurobi_mip_b}
\begin{tabular}{lcccccccc}
\toprule
{Scenario} & \multicolumn{3}{c}{Quantiles} & {Mean} & {Geomean} & {No Sol} & {Wins} \\
{} & {0.1} & {0.5} & {0.9} & {} & {} & {} & {} \\
\midrule
PNFT-0.2-5 & 0.00 & 1.28 & 23.27 & 7.72 & 3.71 & 1 & 9 \\
PNFT-0.0-5 & 0.00 & 2.05 & 22.12 & 7.24 & 3.70 & 0 & 9 \\
PNF-0.2-5 & 0.00 & 1.03 & 17.70 & 6.27 & 3.04 & 1 & 12 \\
PNF-0.5-5 & 0.00 & 1.28 & 29.65 & 8.98 & 3.79 & 1 & 10 \\
RINS-0.2-5 & 0.00 & 4.14 & 27.56 & 14.17 & 4.72 & 1 & 8 \\
RINS-0.5-5 & 0.00 & 6.25 & 37.45 & 13.03 & 5.39 & 2 & 9 \\
LB-0.7-5 & 0.00 & 1.21 & 27.76 & 8.44 & 3.73 & 1 & 12 \\
LB-0.8-5 & 0.00 & 0.24 & 23.93 & 6.60 & 3.06 & 1 & 15 \\
grb-f1 & 0.00 & 0.46 & 15.81 & 6.24 & 3.14 & 1 & 13 \\
grb & 0.00 & 2.97 & 31.03 & 14.17 & 4.48 & 1 & 13 \\
\bottomrule
\end{tabular}
\end{table}

\begin{table}[h]
\centering
\caption{Primal gap (\%) using SCIP with a 20-minute limit on MIPLIB-S1-B}
\label{tab:gap_scip_mip_b}
\begin{tabular}{lcccccccc}
\toprule
{Scenario} & \multicolumn{3}{c}{Quantiles} & {Mean} & {Geomean} & {No Sol} & {Wins} \\
{} & {0.1} & {0.5} & {0.9} & {} & {} & {} & {} \\
\midrule
PNFT-0.2-5 & 0.00 & 19.57 & 87.23 & 37.23 & 12.21 & 8 & 7 \\
PNFT-0.0-5 & 0.00 & 19.57 & 87.23 & 36.89 & 11.84 & 8 & 7 \\
PNF-0.2-5 & 0.00 & 19.62 & 116.28 & 38.05 & 12.34 & 7 & 13 \\
PNF-0.5-5 & 0.00 & 20.59 & 90.54 & 41.51 & 13.44 & 8 & 10 \\
RINS-0.2-5 & 0.00 & 36.79 & 146.03 & 62.97 & 15.04 & 7 & 11 \\
RINS-0.5-5 & 0.00 & 32.50 & 93.77 & 43.22 & 14.10 & 8 & 7 \\
LB-0.7-5 & 0.12 & 31.05 & 111.80 & 176.34 & 19.12 & 5 & 9 \\
LB-0.8-5 & 0.10 & 30.85 & 199.21 & 107.14 & 19.43 & 5 & 8 \\
scip & 0.00 & 14.19 & 98.75 & 111.20 & 14.70 & 2 & 15 \\
\bottomrule
\end{tabular}
\end{table}

\begin{table}[t]
\centering
\caption{Primal gap (\%) and runtime using CPLEX with a 20-minute limit on MIPLIB-S1-B}
\label{tab:gap_time_cplex_miplib_b}
\begin{tabular}{lccccccc}
\toprule
{Instance} & \multicolumn{2}{c}{PNFT-0.2-5} & \multicolumn{2}{c}{RINS-0.2-5} & \multicolumn{2}{c}{cplex} \\
{} & {PG (\%)} & {Time (s)} & {PG (\%)} & {Time (s)} & {PG (\%)} & {Time (s)} \\
\midrule
ns1856153 & \underline{\textbf{99.14}} & 307.73 & \underline{\textbf{99.14}} & 301.29 & \underline{\textbf{99.14}} & \underline{\textbf{138.37}} \\
neos-3634244-kauru & 112.6 & \underline{\textbf{843.28}} & 104.98 & 1135.43 & \underline{\textbf{60.49}} & 1196.69 \\
genus-g31-8 & \underline{\textbf{57.14}} & \underline{\textbf{319.37}} & nan & nan & nan & nan \\
neos-4358725-tagus & \underline{\textbf{17.4}} & 310.77 & 1738.52 & 319.86 & 25.71 & \underline{\textbf{19.12}} \\
neos-3009394-lami & \underline{\textbf{0.0}} & 308.32 & 0.0 & 301.02 & 0.0 & \underline{\textbf{21.52}} \\
ds & \underline{\textbf{28.87}} & \underline{\textbf{469.81}} & 38.08 & 658.78 & 49.48 & 651.74 \\
momentum3 & 120.86 & \underline{\textbf{505.58}} & 103.76 & 697.95 & \underline{\textbf{100.67}} & 658.34 \\
milo-v13-4-3d-4-0 & 2.74 & 1138.4 & \underline{\textbf{1.6}} & \underline{\textbf{301.02}} & 17.55 & 1097.39 \\
neos-3209462-rhin & 66.67 & 742.18 & 66.67 & 748.36 & \underline{\textbf{50.0}} & \underline{\textbf{431.41}} \\
milo-v13-4-3d-3-0 & 7.7 & 1182.58 & \underline{\textbf{0.0}} & \underline{\textbf{346.82}} & 10.14 & 1147.57 \\
genus-sym-grafo5708-48 & 63.16 & 461.65 & \underline{\textbf{42.11}} & 615.78 & 94.74 & \underline{\textbf{1.3}} \\
genus-sym-g31-8 & \underline{\textbf{47.62}} & \underline{\textbf{644.21}} & nan & nan & 57.14 & 1034.01 \\
tw-myciel4 & \underline{\textbf{0.0}} & 307.65 & \underline{\textbf{0.0}} & 301.13 & \underline{\textbf{0.0}} & \underline{\textbf{12.23}} \\
shipsched & 21.73 & 1204.67 & \underline{\textbf{19.2}} & \underline{\textbf{1161.96}} & 30.42 & 1171.35 \\
graphdraw-mainerd & 62.32 & 1196.72 & 89.52 & \underline{\textbf{487.49}} & \underline{\textbf{54.99}} & 1193.03 \\
t1717 & \underline{\textbf{7.79}} & 307.67 & \underline{\textbf{7.79}} & 301.55 & \underline{\textbf{7.79}} & \underline{\textbf{1.5}} \\
neos-5260764-orauea & \underline{\textbf{1.69}} & 1195.74 & 8.65 & 1194.14 & 5.29 & \underline{\textbf{1149.9}} \\
graph20-20-1rand & \underline{\textbf{0.0}} & 307.91 & \underline{\textbf{0.0}} & 301.05 & \underline{\textbf{0.0}} & \underline{\textbf{187.79}} \\
tokyometro & \underline{\textbf{0.0}} & 742.6 & 13.17 & \underline{\textbf{368.5}} & 9.12 & 1196.24 \\
neos-2629914-sudost & 0.09 & 1197.48 & \underline{\textbf{0.01}} & 1174.79 & 0.15 & \underline{\textbf{1133.37}} \\
neos-1067731 & 0.0 & 1147.4 & 0.02 & 1114.91 & \underline{\textbf{0.0}} & \underline{\textbf{1088.34}} \\
bley\_xs1noM & \underline{\textbf{0.0}} & 1175.9 & 8.3 & \underline{\textbf{974.78}} & 3.18 & 1177.16 \\
assign1-10-4 & \underline{\textbf{0.24}} & 307.46 & \underline{\textbf{0.24}} & 301.11 & \underline{\textbf{0.24}} & \underline{\textbf{15.26}} \\
t1722 & \underline{\textbf{6.98}} & 1176.76 & 9.49 & \underline{\textbf{301.32}} & 9.46 & 1150.14 \\
usafa & nan & nan & nan & nan & \underline{\textbf{17446.43}} & \underline{\textbf{1.47}} \\
neos-2294525-abba & nan & nan & nan & nan & nan & nan \\
neos-3581454-haast & \underline{\textbf{0.0}} & \underline{\textbf{456.88}} & 0.74 & 469.2 & \underline{\textbf{0.0}} & 515.65 \\
bppc6-06 & 3.29 & \underline{\textbf{317.26}} & \underline{\textbf{0.94}} & 423.96 & \underline{\textbf{0.94}} & 1154.14 \\
eilA101-2 & \underline{\textbf{2.82}} & 325.86 & 3.52 & \underline{\textbf{301.83}} & 3.37 & 400.28 \\
neos-631710 & \underline{\textbf{0.0}} & 309.11 & \underline{\textbf{0.0}} & 302.8 & \underline{\textbf{0.0}} & \underline{\textbf{90.76}} \\
\bottomrule
\end{tabular}
\end{table}

\begin{table}[t]
\centering
\caption{Primal gap and runtime using Gurobi with a 20-minute limit on MIPLIB-S1-B}
\label{tab:primal_gap_time_grb_miplib_b}
\begin{tabular}{lccccccc}
\toprule
{Instance} & \multicolumn{2}{c}{PNFT-0.2-5} & \multicolumn{2}{c}{RINS-0.2-5} & \multicolumn{2}{c}{grb} \\
{} & {PG (\%)} & {Time (s)} & {PG (\%)} & {Time (s)} & {PG (\%)} & {Time (s)} \\
\midrule
ns1856153 & \underline{\textbf{0.0}} & 901.42 & 48.44 & \underline{\textbf{300.41}} & 24.27 & 757.2 \\
neos-3634244-kauru & 16.01 & 1108.62 & 9.93 & 817.56 & \underline{\textbf{0.0}} & \underline{\textbf{747.51}} \\
genus-g31-8 & \underline{\textbf{9.52}} & 300.08 & \underline{\textbf{9.52}} & 300.08 & 28.57 & \underline{\textbf{168.02}} \\
neos-4358725-tagus & \underline{\textbf{0.32}} & 707.19 & 195.8 & 302.09 & 2.53 & \underline{\textbf{153.14}} \\
neos-3009394-lami & \underline{\textbf{0.0}} & 332.15 & 0.0 & 300.36 & \underline{\textbf{0.0}} & \underline{\textbf{15.45}} \\
ds & \underline{\textbf{24.46}} & 390.37 & 32.43 & 354.99 & 30.15 & \underline{\textbf{275.85}} \\
momentum3 & 34.6 & 945.4 & \underline{\textbf{23.52}} & \underline{\textbf{899.18}} & 47.39 & 1143.05 \\
milo-v13-4-3d-4-0 & 4.14 & 300.38 & 4.14 & 300.41 & \underline{\textbf{0.0}} & \underline{\textbf{125.71}} \\
neos-3209462-rhin & \underline{\textbf{16.67}} & \underline{\textbf{300.23}} & \underline{\textbf{16.67}} & \underline{\textbf{300.23}} & \underline{\textbf{16.67}} & 494.28 \\
milo-v13-4-3d-3-0 & 11.63 & 391.36 & 15.22 & 300.37 & \underline{\textbf{3.84}} & \underline{\textbf{228.04}} \\
genus-sym-grafo5708-48 & 42.11 & 664.13 & \underline{\textbf{0.0}} & \underline{\textbf{443.08}} & 21.05 & 572.0 \\
genus-sym-g31-8 & 38.1 & \underline{\textbf{300.09}} & \underline{\textbf{28.57}} & \underline{\textbf{300.08}} & 57.14 & 320.0 \\
tw-myciel4 & \underline{\textbf{0.0}} & 300.07 & \underline{\textbf{0.0}} & 300.07 & \underline{\textbf{0.0}} & \underline{\textbf{0.61}} \\
shipsched & 9.18 & \underline{\textbf{300.22}} & 9.17 & 367.7 & \underline{\textbf{0.0}} & 1196.02 \\
graphdraw-mainerd & \underline{\textbf{18.54}} & \underline{\textbf{753.12}} & 21.35 & 887.0 & 31.25 & 1198.79 \\
t1717 & \underline{\textbf{8.36}} & 300.09 & \underline{\textbf{8.36}} & 300.09 & \underline{\textbf{8.36}} & \underline{\textbf{0.15}} \\
neos-5260764-orauea & 0.46 & \underline{\textbf{300.17}} & 0.46 & 1087.41 & \underline{\textbf{0.0}} & 421.81 \\
graph20-20-1rand & \underline{\textbf{0.0}} & 300.08 & \underline{\textbf{0.0}} & 300.08 & \underline{\textbf{0.0}} & \underline{\textbf{2.41}} \\
tokyometro & \underline{\textbf{1.06}} & 431.47 & 16.46 & \underline{\textbf{300.07}} & 2.99 & 507.03 \\
neos-2629914-sudost & \underline{\textbf{0.02}} & 1026.36 & 0.03 & 1146.71 & 0.04 & \underline{\textbf{943.26}} \\
neos-1067731 & \underline{\textbf{0.0}} & \underline{\textbf{493.84}} & 0.02 & 1008.18 & 0.03 & 963.72 \\
bley\_xs1noM & 3.42 & 1136.62 & 7.0 & \underline{\textbf{1059.12}} & \underline{\textbf{2.86}} & 1123.1 \\
assign1-10-4 & 0.24 & 300.09 & 0.24 & 300.36 & \underline{\textbf{0.0}} & \underline{\textbf{45.98}} \\
t1722 & 1.28 & 300.08 & \underline{\textbf{1.1}} & 300.08 & 7.67 & \underline{\textbf{274.55}} \\
usafa & \underline{\textbf{0.55}} & 1064.7 & 4.57 & \underline{\textbf{633.0}} & 5.18 & 1094.74 \\
neos-2294525-abba & nan & nan & nan & nan & nan & nan \\
neos-3581454-haast & \underline{\textbf{0.0}} & 300.13 & \underline{\textbf{0.0}} & 300.38 & \underline{\textbf{0.0}} & \underline{\textbf{20.02}} \\
bppc6-06 & \underline{\textbf{0.0}} & 300.16 & \underline{\textbf{0.0}} & 300.39 & 3.29 & \underline{\textbf{33.01}} \\
eilA101-2 & 2.82 & 300.09 & \underline{\textbf{0.25}} & 312.33 & 2.97 & \underline{\textbf{76.03}} \\
neos-631710 & \underline{\textbf{0.0}} & 300.1 & \underline{\textbf{0.0}} & 300.1 & \underline{\textbf{0.0}} & \underline{\textbf{82.83}} \\
\bottomrule
\end{tabular}
\end{table}

\begin{table}[t]
\centering
\caption{Primal integral and solution count using Gurobi with a 20-minute limit on MIPLIB-S1-B}
\label{tab:primal_integral}
\begin{tabular}{lccccccc}
\toprule
{Instance} & \multicolumn{2}{c}{PNFT-0.2-5} & \multicolumn{2}{c}{RINS-0.2-5} & \multicolumn{2}{c}{grb} \\
{} & {PI} & {Sol count} & {PI} & {Sol count} & {PI} & {Sol count} \\
\midrule
ns1856153 & 13.88 & 6 & \underline{\textbf{0.05}} & \underline{\textbf{3}} & 43.91 & 6 \\
neos-3634244-kauru & 24.92 & 13 & \underline{\textbf{23.19}} & 18 & 24.54 & \underline{\textbf{10}} \\
genus-g31-8 & \underline{\textbf{0.0}} & \underline{\textbf{1}} & \underline{\textbf{0.0}} & \underline{\textbf{1}} & 45.2 & 5 \\
neos-4358725-tagus & 95.42 & 22 & 150.81 & \underline{\textbf{16}} & \underline{\textbf{81.87}} & 41 \\
neos-3009394-lami & \underline{\textbf{0.0}} & 4 & 0.0 & \underline{\textbf{3}} & 0.0 & 8 \\
ds & \underline{\textbf{19.51}} & 5 & 30.52 & \underline{\textbf{3}} & 149.87 & 11 \\
momentum3 & 30.38 & \underline{\textbf{14}} & \underline{\textbf{25.61}} & 19 & 66.94 & \underline{\textbf{14}} \\
milo-v13-4-3d-4-0 & 3.11 & \underline{\textbf{3}} & 3.11 & \underline{\textbf{3}} & \underline{\textbf{1.21}} & 12 \\
neos-3209462-rhin & \underline{\textbf{0.02}} & \underline{\textbf{2}} & 0.02 & \underline{\textbf{2}} & 62.81 & 6 \\
milo-v13-4-3d-3-0 & 8.79 & 5 & 8.91 & \underline{\textbf{3}} & \underline{\textbf{4.52}} & 8 \\
genus-sym-grafo5708-48 & 34.42 & \underline{\textbf{3}} & \underline{\textbf{17.09}} & 6 & 57.34 & 8 \\
genus-sym-g31-8 & \underline{\textbf{0.0}} & \underline{\textbf{1}} & 18.78 & 2 & 79.11 & 5 \\
tw-myciel4 & \underline{\textbf{0.0}} & \underline{\textbf{1}} & \underline{\textbf{0.0}} & \underline{\textbf{1}} & 13.93 & 3 \\
shipsched & \underline{\textbf{0.01}} & \underline{\textbf{2}} & 6.88 & 4 & 10.39 & 48 \\
graphdraw-mainerd & \underline{\textbf{15.53}} & 28 & 17.5 & \underline{\textbf{25}} & 66.6 & 59 \\
t1717 & \underline{\textbf{0.0}} & \underline{\textbf{1}} & \underline{\textbf{0.0}} & \underline{\textbf{1}} & \underline{\textbf{0.0}} & \underline{\textbf{1}} \\
neos-5260764-orauea & \underline{\textbf{0.0}} & \underline{\textbf{2}} & 0.34 & 5 & 14.15 & 27 \\
graph20-20-1rand & \underline{\textbf{0.0}} & \underline{\textbf{1}} & \underline{\textbf{0.0}} & \underline{\textbf{1}} & 5.31 & 5 \\
tokyometro & 1.57 & 18 & \underline{\textbf{0.12}} & \underline{\textbf{2}} & 41.82 & 53 \\
neos-2629914-sudost & \underline{\textbf{0.11}} & \underline{\textbf{8}} & 0.15 & 10 & 0.28 & 12 \\
neos-1067731 & \underline{\textbf{0.0}} & 21 & 0.01 & \underline{\textbf{17}} & 9.36 & 65 \\
bley\_xs1noM & \underline{\textbf{4.75}} & 34 & 5.72 & \underline{\textbf{22}} & 10.19 & 108 \\
assign1-10-4 & \underline{\textbf{0.0}} & \underline{\textbf{2}} & 0.18 & 3 & 1.72 & 15 \\
t1722 & \underline{\textbf{0.0}} & \underline{\textbf{1}} & 0.5 & 2 & 9.76 & 30 \\
usafa & 6.25 & \underline{\textbf{15}} & \underline{\textbf{5.82}} & \underline{\textbf{15}} & 1106.55 & 17 \\
neos-2294525-abba & nan & \underline{\textbf{0}} & nan & \underline{\textbf{0}} & nan & \underline{\textbf{0}} \\
neos-3581454-haast & \underline{\textbf{0.0}} & \underline{\textbf{2}} & \underline{\textbf{0.0}} & 3 & 5.43 & 12 \\
bppc6-06 & \underline{\textbf{0.0}} & \underline{\textbf{2}} & \underline{\textbf{0.0}} & 3 & 3.54 & 13 \\
eilA101-2 & \underline{\textbf{0.0}} & \underline{\textbf{1}} & 0.33 & 3 & 33.01 & 5 \\
neos-631710 & \underline{\textbf{0.0}} & \underline{\textbf{1}} & \underline{\textbf{0.0}} & \underline{\textbf{1}} & 67.61 & 13 \\
\bottomrule
\end{tabular}
\end{table}

\newpage

\clearpage

\newpage

\bibliography{mybib}

\bibliographystyle{plainnat}

\end{document}